\documentclass[12pt,A4]{article}
\usepackage{theorem}
\newtheorem{lemma}{Lemma}
\newtheorem{proposition}[lemma]{Proposition}
\newtheorem{theorem}[lemma]{Theorem}
\newtheorem{corollary}[lemma]{Corollary}
{\theorembodyfont{\upshape}}
{\theorembodyfont{\upshape}}
{\theorembodyfont{\upshape}}
{\theorembodyfont{\upshape}}
{\theorembodyfont{\upshape}}

\newcommand{\Z}{{\bf Z}}
\newcommand{\R}{{\bf R}}
\newcommand{\C}{{\bf C}}
\newcommand{\rme}{{\rm e}}

\newcommand{\cH}{{\cal H}}
\newcommand{\cK}{{\cal K}}
\newcommand{\cL}{{\cal L}}
\newcommand{\cM}{{\cal M}}
\newcommand{\cN}{{\cal N}}

\newcommand{\cQ}{{\cal Q}}

\newcommand{\sig}{\sigma}
\newcommand{\alp}{\alpha}
\newcommand{\bet}{\beta}
\newcommand{\gam}{\gamma}

\newcommand{\lam}{\lambda}
\newcommand{\del}{\delta}
\newcommand{\eps}{\varepsilon}

\newcommand{\Gam}{\Gamma}
\newcommand{\Ome}{\Omega}

\newcommand{\Spec}{{\rm Spec}}

\newcommand{\norm}{\Vert}

\newcommand{\Proof}{\underbar{Proof}{\hskip 0.1in}}

\newcommand{\cNum}{{\rm Num}}
\newcommand{\conv}{{\rm Conv}}
\newcommand{\cconv}{{\rm Conv}}
\newcommand{\dist}{{\rm dist}}
\newcommand{\Schrodinger}{Schr\"odinger }

\newcommand{\sigl}{\sig_{{\rm loc}}}
\newcommand{\var}{{\rm var}}
\newcommand{\tensor}{\otimes}
\newcommand{\inv}{^{-1}}
\newcommand{\half}{\frac{1}{2}}
\newcommand{\mod}{{\rm \,\,mod\,\,}}
\title{SPECTRAL THEORY\\ OF PSEUDO-ERGODIC OPERATORS}
\author{E.B. Davies}
\date{June 2000}
\begin{document}
\maketitle
\begin{abstract}
We define a class of pseudo-ergodic non-self-adjoint
\Schrodinger operators acting in spaces $l^2(X)$ and
prove some general theorems about their spectral properties. We
then apply these to study the spectrum of a non-self-adjoint
Anderson model acting on $l^2(\Z)$, and find the precise condition
for $0$ to lie in the spectrum of the operator. We also introduce
the notion of localized spectrum for such operators.
\par
keywords: Anderson model, spectrum, random, ergodic
\Schrodinger operator, pseudospectrum, non-self-adjoint operator.
\end{abstract}
\section{Introduction}
\par
Recent papers have obtained some striking results concerning the
spectral properties of the non-self-adjoint (nsa) Anderson model,
which models the growth of bacteria in an inhomogeneous
environment, \cite{Gold,Hatano1,Hatano2,Nelson,Dahmen}. To be more
precise the authors have determined the asymptotic limit of the
spectrum of a nsa random finite periodic chain almost surely as the
length of the chain increases to infinity. In a later paper the
author considered the same random operator $H$ acting on $l^2(\Z)$,
and found that the spectrum is very different from that obtained by
the cited authors, \cite{EBD}. The reason for this is that the
spectral properties of nsa operators are highly unstable, and
infinite volume limits should be examined using pseudospectral
ideas, \cite{Boe,Boe2,BS,EBD,Dav4,Dav3,Red,RT,Tre1,Tre2}. More
specifically if $\lam$ lies in the spectrum of the infinite volume
nsa Anderson model, it need not be close to the spectrum of the
finite volume periodic Anderson model; one expects rather that the
norm of the resolvent operator $(H-\lam )\inv$ of the finite volume
model will diverge as the volume increases. These pseudospectral
ideas have been worked out in detail for a random {\it bidiagonal}
model, which is in a certain sense exactly soluble,
\cite{BZ,FZ,TCE}. Our results may therefore be interpreted as
finding the region in the complex plane for which the finite volume
nsa periodic Anderson model has very large resolvent norm.

In the present paper we reconsider such problems in a more general
context, in which the probabilistic aspects have been eliminated in
favour of what we call pseudo-ergodic ideas. As well as making the
subject more accessible to those without a probabilistic training,
this emphasizes the fact that the spectral matters which we
consider depend only on the support of the relevant probability
measure. On the other hand the asymptotics of the spectrum of the
finite volume periodic nsa Anderson model does depend on the
probability measure. We finally carry out a more detailed spectral
analysis of the infinite volume nsa Anderson operator, and find
precise conditions under which zero almost surely lies in the
spectrum. We also obtain further results on the location of the
spectrum, which come close to a complete determination in many
cases. In the final section we consider the possibility that there
may be constraints on the pair of values of the potential at two
neighbouring points which are absolute rather than just
probabilistic.

\section{The general context}

The operators which we consider act on the Hilbert space
$l^2(X,\cK)\sim l^2(X)\tensor \cK$, where $X$ is a countable set on
which a group $\Gam$ acts by permutations. The simplest choice of
the auxiliary Hilbert space $\cK$ is $\C$, but other choices are
needed in some applications; see the end of Section 3. Many of the
results presented here apply to $l^p(X,\cK)$ with $p\not= 2$
without modification (the case $p=1$ is of probabilistic
importance), but this does not apply to those involving numerical
ranges. We define the unitary operators $U_\gam$ for $\gam\in\Gam$
by $U_\gam f(x)=f(\gam^{-1}x)$. The bounded operators which we
study are of the form $H=H_0\tensor I+V$. Here $H_0$ acts on
$l^2(X)$ and commutes with the action of $\Gam$ in the sense that
$H_0 U_\gam
=U_\gam H_0$ for all $\gam\in\Gam$, or equivalently
\[
H_0(\gam x,\gam y)=H_0(x,y)
\]
for all $\gam\in\Gam$ and all $x,y\in X$, where $H_0(x,y)$ is the
infinite matrix associated with $H_0$. We assume that the spectrum
$E$ of $H_0$ is known. From this point onwards we write $H_0$ for
$H_0\tensor I$.

Given a norm closed, bounded set $\cM\subseteq \cL(\cK)$, we assume
that the operator $V$ is of the form
\[
(Vf)(x)=V(x)f(x)
\]
where $V(x)\in\cM$ for all $x\in X$. We say that $V$ is
$(\Gam,\cM)$ pseudo-ergodic if its set of spatial translates is
dense in the following sense. For every $\eps
>0$, every finite subset $F\subset X$ and every $W:F\to \cM$, there
exists $\gam\in\Gam$ such that
\[
\norm W(x)-V(\gam x)\norm <\eps
\]
for all $x\in F$. It is well known that a large class of suitably
defined random potentials have this property almost surely, but we
consider a single potential, and do not need to introduce any
probabilistic ideas. The same class of pseudo-ergodic potentials is
applicable to a variety of different random models, as we explain
in more detail in the final section.

The above definition suffices for our purposes, but it does not
capture the full sense of random behaviour and may be refined as
follows. We define a direction $U$ to be an infinite subset of $X$
such that for every finite $F\subset X$ there exists $\gam\in\Gam$
such that $\gam F\subset U$. We then say that $V$ is $(\Gam,\cM)$
pseudo-ergodic in the direction $U$ if for every $\eps
>0$, every finite subset $F\subset X$ and every $W:F\to \cM$, there
exists $\gam\in\Gam$ such that $\gam F\subset U$ and
\[
\norm W(x)-V(\gam x)\norm <\eps
\]
for all $x\in F$. Suitably defined random potentials have this
property for every choice of direction almost surely, and therefore
have the property simultaneously for any countable set of
directions almost surely. The property itself, however, is defined
for a single potential and makes no mention of probability.

The following theorem is an adaptation of a well-known result of
Pastur for random potentials. We will use it to approximate
$\Spec(H)$ from inside by making suitable choices of $W$.

\begin{theorem} If $H=H_0+V$ where $V$ is $(\Gam,\cM)$ pseudo-ergodic
and $K=H_0+W$ where $W:X\to \cM$ is arbitrary, then
\[
\Spec(K)\subseteq \Spec(H).
\]
In particular if $V,W$ are both $(\Gam,\cM)$ pseudo-ergodic then
they have the same spectrum.
\end{theorem}

\Proof If $\lam\in\Spec(K)$ then there exists a sequence $f_n\in
l^2(X,\cK)$ with $\norm f_n \norm =1$ and either $\norm Kf_n-\lam
f_n\norm\to 0$ or $\norm K^\ast f_n-\overline{\lam} f_n\norm\to 0$;
we consider only the former case, the latter being similar. Given
$\eps >0$ a truncation procedure shows that there exists $f$ with
finite support $F$ in $X$ such that $\norm f\norm=1$ and $\norm K
f-\lam f\norm <\eps/2$. Since $V$ is pseudo-ergodic there exists
$\gam\in\Gam$ such that $\norm H_\gam f-Kf\norm<\eps/2$, where
\[
H_\gam=U_\gam^{-1}HU_\gam=H_0+V(\gam\,\cdot).
\]
Putting $f_\eps=U_\gam f$ we deduce that
\[
\norm H f_\eps -\lam f_\eps\norm =\norm U_\gam^{-1} H U_\gam
f-\lam f\norm <\eps
\]
and the arbitrariness of $\eps >0$ implies that $\lam\in\Spec(H)$.

\begin{corollary} If $H=H_0+V$ where $V$ is $(\Gam,\cM)$ pseudo-ergodic
then
\[
\Spec(H)=\bigcup \{ \Spec(H_0+W):W\in \cM^X\}.
\]
If also $\tilde H=H_0+\tilde V$ where $\tilde V$ is $(\Gam,\tilde
\cM)$ pseudo-ergodic with $\cM\subseteq \tilde \cM$ then
\[
\Spec(H)\subseteq \Spec(\tilde H).
\]
\end{corollary}

From this point we assume that $H=H_0+V$ where $V$ is $(\Gam,\cM)$
pseudo-ergodic. We put
\[
\Spec(\cM)=\overline{\bigcup_{A\in\cM}\Spec(A)}
\]
and
\[
\cNum(\cM)=\overline{\bigcup_{A\in\cM}\cNum(A)}
\]
where $\cNum$ denotes the closure of the numerical range.

\begin{theorem}
The spectrum of $H$ satisfies
\[
E+\Spec(\cM) \subseteq
\Spec(H)\subseteq
\cNum(H_0)+\cconv(\cNum(\cM))
\]
where $\cconv$ denotes the closed convex hull. If $H_0$ is normal
and $A$ is normal for every $A\in\cM$ then
\begin{equation}
\Spec(H)\subseteq \cconv(E)+\cconv(\Spec (\cM))\label{locatespec}
\end{equation}
\end{theorem}

\Proof Theorem 1 implies that for each $A\in \cM$
\[
E+\Spec(A)=\Spec(H_0\tensor I+I\tensor A)\subseteq \Spec(H)
\]
and this yields the first inclusion. The second depends on use of
the numerical range to give
\begin{eqnarray*}
\Spec(H)&\subseteq &\cNum(H)\\
&\subseteq& \cNum(H_0)+\cNum(V).
\end{eqnarray*}
Now $z$ lies in the numerical range of $V$ if and only if there
exists $f\in l^2(X,\cK)$ of norm $1$ such that $z=\langle
Vf,f\rangle$. Putting $g_x=f(x)/\norm f(x)\norm$, provided this is
non-zero, and $\mu_x=\norm f(x)\norm^2$, we see that $\mu$ is a
probability measure on $X$ and that
\[
z=\sum_{x\in X} \mu_x \langle V_xg_x,g_x\rangle\in
\cconv(\cNum(M)).
\]
Hence $\cNum(V)\subseteq \cconv(\cNum(M))$, and the first statement
of the theorem follows. The second statement is a consequence of
the fact that $\cNum(B)$ equals $\cconv(\Spec(B))$ for any normal
operator $B$.

Let $B(x,r)$ denote the closed ball $\{ y:|x-y|\leq r\}$. The next
theorem complements Theorem 3.

\begin{theorem}
If $A$ is normal for every $A\in\cM$ then the spectrum of $H$
satisfies
\begin{equation}
\Spec(H)\subseteq \Spec(\cM)+B(0,e)\label{4a}
\end{equation}
where $e=\norm H_0 \norm$. If $H_0$ is normal then
\begin{equation}
\Spec(H)\subseteq E+B(0,\mu)\label{4b}
\end{equation}
where $\mu=max\{\norm A \norm: A\in \cM\}$.
\end{theorem}

\Proof If $V$ is normal then using $\Spec (V)\subseteq \Spec(\cM)$
we see that
\[
\norm (V-z I)\inv\norm =\dist (z,\Spec(V))\inv\leq \dist (z,\Spec(\cM))\inv
\]
for all $z\notin \cM$. Since $z\notin \cM+B(0,e)$ is equivalent to
$\dist\{ z,\Spec(\cM)\}>\norm H_0\norm$, it implies
\[
\norm H_0 \norm\norm (V-z I)\inv\norm <1
\]
and the resolvent expansion for $(H_0+V-z I)\inv$ is norm
convergent. The proof of the second part of the theorem is similar.

We also wish to classify the spectrum of nsa operators acting on
$l^2(X,\cK)$, and for this purpose we assume that $X$ is provided
with a metric $d$ such that every ball $B(x,r)=\{y\in X :d(x,y)\leq
r)\}$ is finite and such that $\Gam$ acts as a group of isometries
of $X$. Given a function $f:X\to \cK$ with $\norm f\norm_2=1$ we
define its variance by
\[
\var(f)=\min_{y\in X}\sum_{x\in X}d(x,y)^2 |f(x)|^2
\]
and its expectation to be any of the points in $X$ at which the
minimum is achieved. The following theorems have analogues in which
the variance is replaced by higher order moments, or suitable
subexponential weights.

\begin{lemma}
If $\norm f\norm_2=1$ and
\[
v(x)=\sum_{y\in X}d(x,y)^2 |f(y)|^2
\]
is finite for some $x\in X$ then it is finite for every $x\in X$
and $v(x)$ increases indefinitely as $x\to\infty$. Thus the minimum
of $v(\cdot )$ is achieved at a finite number of points only. If
$x_i$, $i=1,2$, are points at which $v$ has the same minimum value
$s$ then $d(x_1,x_2)\leq 2s^{1/2}$.
\end{lemma}

\Proof If $v(x)<\infty$ then for any $u\in X$ we have
\[
v(u)\leq 2\sum_{y\in
X}\{d(x,y)^2+d(x,u)^2\}|f(y)|^2=2\{v(x)+d(x,u)^2\} <\infty.
\]
by the triangle inequality. If the finite set $F$ satisfies
\[
\sum_{y\in F}|f(y)|^2\geq \frac{1}{2}
\]
then
\[
v(x)\geq \sum_{y\in F} d(x,y)^2 |f(y)|^2\geq \frac{1}{2} d(x,F)^2
\]
which increases indefinitely as $x\to\infty$ because of our
assumption that all balls of finite radius contain only a finite
number of points.

Now suppose that $s=\min\{v(x):x\in X\}$ and that
$v(x_1)=v(x_2)=s$. Then by the triangle inequality
\begin{eqnarray*}
2s&=&\sum_{y\in X}\{ d(x_1,y)^2+d(x_2,y)^2\} |f(y)|^2\\
&\geq&\half\sum_{y\in X} d(x_1,x_2)^2 |f(y)|^2\\ &=&\half
d(x_1,x_2)^2
\end{eqnarray*}
which implies the second statement of the lemma.

Following \cite{EBD} we define the localized spectrum $\sigl (A)$
of any bounded operator $A$ on $l^2(X,\cK)$ to be the set of all
$\lam\in
\C$ such that there exists a sequence $f_n\in l^2(X,\cK)$ of unit
vectors such that $\norm Af_n-\lam f_n\norm\to 0$ while $\var
(f_n)$ remains uniformly bounded. If $\lam$ is an eigenvalue then
one would expect its corresponding eigenfunction to decrease
rapidly at infinity and hence to have finite variance, in which
case $\lam$ would lie in $\sigl (A)$. What is more surprising is
that $\sigl (A)$ can be much larger than the set of eigenvalues of
$A$.

\begin{theorem}
If $H=H_0+V$ where $V$ is $(\Gam,\cM)$ pseudo-ergodic and $K=H_0+W$
where $W:X\to \cM$ is arbitrary, then
\[
\sigl(K)\subseteq \sigl(H).
\]
Thus every eigenvalue of $K$ lies in the localized spectrum of $H$.
Moreover if $V,W$ are both $(\Gam,\cM)$ pseudo-ergodic then they
have the same localized spectrum.
\end{theorem}

\Proof First note that if $f\in l^2(X,\cK)$ has unit norm and $\gam\in\Gam$ then
$g=U_\gam f$ has the same variance as $f$ because $\Gam$ acts as a
group of isometries of $X$. It is a consequence of the definition
of pseudo-ergodicity that there exists a sequence $\gam(n)\in\Gam$
such that $H_n=U_{\gam(n)}\inv HU_{\gam(n)}$ converges strongly to
$K$. Now let $\norm f_m\norm =1$, $\var(f_m)\leq s$ and $\norm K
f_m-\lam f_m\norm <\frac{1}{m}$ for all $m\in \Z^+$. Given $m$
\begin{eqnarray*}
\norm H(U_{\gam(n)} f_m)-\lam(U_{\gam(n)}f_m)\norm
&=& \norm U_{\gam(n)}\inv HU_{\gam(n)} f_m-\lam f_m \norm\\
&=&\norm H_n f_m-\lam f_m\norm \\ &\to& \norm K f_m-\lam f_m \norm
<\frac{1}{m}
\end{eqnarray*}
as $n\to\infty$. Therefore there exists $n(m)$ such that
$g_m=U_{\gam(n(m))}f_m$ satisfies
\[
\norm Hg_m-\lam g_m \norm<\frac{1}{m}\]
for all $m\in \Z^+$. Since $\var(g_m)\leq s$ for all $m$ it follows
that $\lam\in\sigl (H)$.

We next turn to the essential spectrum. We say that $z$ lies in the
essential spectrum of a bounded operator $A$ if $A-zI$ is not a
Fredholm operator. We will need the following known result.

\begin{proposition}
Suppose that $z\in\C$ and for all $\eps >0$ and all finite $N$
there exists an orthonormal set $f_1,...,f_N$ such that $\norm
Af_n-zf_n\norm <\eps$ for all $1\leq n\leq N$. Then $z$ lies in the
essential spectrum of $A$.
\end{proposition}

\Proof
Suppose that $z\in\C$ satisfies the conditions of the proposition.
If $\ker(A-zI)$ is infinite dimensional then $A-zI$ is obviously
not Fredholm, so let $\dim(\ker(A-zI))<N$ where $N$ is finite. The
assumption implies that for all $\eps
>0$ there exists an $N$-dimensional subspace $L$
such that $f\in L$ implies
\begin{equation}
\norm A f-zf\norm <\eps \norm f \norm .\label{small}
\end{equation}
Because $\dim(L)>\dim(\ker(A-zI))$ there exists $f\perp \ker(A-zI)$
such that (\ref{small}) holds. Since $\eps >0$ is arbitrary, $A-zI$
cannot be Fredholm.

\begin{lemma}
Suppose that there exists a $(\Gam,\cM)$ pseudo-ergodic potential
$V$ on $X$ where $\cM\subseteq \cL(\cK)$ contains more than one
point. Then for any finite subset $F$ of $X$ and any finite $N$
there exist $\gam_1,...,\gam_N\in\Gam$ such that
$\{\gam_nF\}_{n=1}^N$ are pairwise disjoint.
\end{lemma}

\Proof Let us first put $N=2$. Let $m_1,m_2\in\cM$ and $\norm
m_1-m_2\norm =2\del>0$. Also let $W:F\to \cL(\cK)$ satisfy
$W(x)=m_i$ for all $x\in F$. Since $V$ is $(\Gam,\cM)$
pseudo-ergodic there exist $\gam_i\in\Gam$ such that $\norm
V(\gam_i x)-m_i\norm <\del$ for all $x\in F$, or equivalently
$\norm V(y)-m_i\norm <\del$ for all $y\in \gam_iF$. This implies
that $\gam_1F\cap
\gam_2F=\emptyset$.

We next prove that if the lemma holds for $N$ then it holds for
$2N$; we can then complete the proof by the use of induction. We
put $\tilde{F}=\bigcup_{j=1}^N \gam_j F$ and let
$\bet_1,\bet_2\in\Gam$ be such that
$\bet_1\tilde{F}\cap\bet_2\tilde{F}=\emptyset$. This yields the
statement of the lemma for the sets $\bet_i\gam_j F$ where $i=1,2$
and $1\leq j\leq N$.

\begin{theorem}
If $H=H_0+V$ where $V$ is $(\Gam,\cM)$ pseudo-ergodic and $\cM$
contains more than one point, then $H$ has no inessential spectrum.
\end{theorem}

\Proof If $\lam\in\Spec(H)$ then either (i) for every $\eps >0$
there exists $f\in l^2(X,\cK)$ such that $\norm f \norm =1$ and
$\norm Hf-\lam f\norm <\eps$, or (ii) for every $\eps >0$ there
exists $f\in l^2(X,\cK)$ such that $\norm f \norm =1$ and $\norm
H^\ast f-\overline\lam f\norm <\eps$. We assume (i), the proof for
(ii) being similar. By approximation we may assume that each $f$
has finite support $F$. Now for any $\eps >0$ and any finite $N$
let $\gam_1,...,\gam_N\in\Gam$ be such that $\gam_i F$ are pairwise
disjoint. Put $\tilde{F}=\bigcup_{i=1}^N \gam_iF$ and define
$W:\tilde{F}\to M$ by $W(\gam_i x)=V(x)$ for all $x\in F$. Since
$V$ is $(\Gam,\cM)$ pseudo-ergodic there exists $\gam\in \Gam$ such
that
\[
\norm V(\gam y)-W(y)\norm <\eps
\]
for all $y\in \tilde{F}$. Thus
\begin{equation}
\norm V(\gam \gam_i x)-V( x)\norm <\eps \label{app}
\end{equation}
for all $x\in F$ and $1\leq i\leq N$.

We now put $f_i(x)=f(\gam_i^{-1}\gam^{-1}x)$ for all $x\in X$ and
observe that $f_i$ have supports within $\gam\gam_i F$, which are
disjoint, so $\{f_i \}_{i=1}^N$ form an orthonormal set. It follows
from condition (i) and (\ref{app}) that
\[
\norm Hf_i-\lam f_i\norm <2\eps
\]
for all $1\leq i\leq N$. This implies that $\lam$ lies in the
essential spectrum of $H$ by Proposition 7.

\section{The nsa Anderson model}

In this section we apply the above ideas to an example of physical
and biological importance. We first consider the one-dimensional
nsa Anderson operator
\begin{equation}
Hf_n=\rme^{-g}f_{n-1}+\rme^gf_{n+1}+V_nf_n\label{schr}
\end{equation}
acting on $l^2(\Z)$ (so that $\cK=\C$), where $g>0$ and $V$ is a
$(\Gam,M)$ pseudo-ergodic potential, $\Gam$ being the group of all
translations of $\Z$ and $M$ being a compact subset of $\C$. The
potential $V$ may be generated by assuming that its values at
different points are independent and identically distributed
according to a probability law which has compact support $M$.

Fourier analysis quickly establishes that $H_0$ is normal with
spectrum the ellipse
\begin{equation}
E=\{\rme^{g+i\theta}+\rme^{-g-i\theta}:\theta\in[0,2\pi]\}\label{ellipse}
\end{equation}
following which Theorem 3 implies that
\begin{equation}
 E+M \subseteq \Spec(H)\subseteq
\conv(E)+\conv (M).\label{nn}
\end{equation}

A more precise determination of $\Spec(H)$ depends upon  the size
of $g$, the choice of $M$ and the use of Theorem 6, extending what
we already proved in \cite{EBD}. Given any finite sequence
$\alp=(\alp_0,\alp_1,\ldots,\alp_{n-1})\in M^n$ let $W_\alp$ be the
periodic potential such that $W_{\alp,m}=\alp_r$ if $m=r\mod n$.
The eigenvalue equation
\begin{equation}
\rme^{-g} f_{m-1} +W_{\alp,m} f_m +\rme^g f_{m+1}=\lam f_m \label{import}
\end{equation}
may be rewritten in terms of $w_m=(f_{m-1},f_m)\in\C^2$ as
$w_{m+1}=w_mA_m$ where
\[
A_m=\left[\begin{array}{cc} 0&-\rme^{-2g}\\
1&\rme^{-g}(\lam-W_{\alp,m})
\end{array} \right].
\]
Thus
\[
w_{n(r+1)}=w_{nr}B
\]
for all $r\in\Z$ where $B$ is the transfer matrix
\[
B=A_0A_1\ldots A_{n-1}.
\]
Since
\[
\det(B)=\prod_{r=0}^{n-1} \det(A_r)=\rme^{-2ng}
\]
it follows that at least one of the two eigenvalues $\mu_1,\mu_2$
of $B$ satisfies $|\mu_i|<1$. If we write
\[
B=\left[\begin{array}{cc} b_{11}(\lam)&b_{12}(\lam)\\
b_{21}(\lam)&b_{22}(\lam)
\end{array} \right]
\]
then one may prove by induction that $b_{22}(\lam)$ is a polynomial
of degree $n$ in $\lam$ while the other coefficients are of lower
degree.

The solution $f$ of (\ref{import}) corresponding to an eigenvalue
$\mu$ of $B$ is exponentially increasing or decreasing on $\Z$
according to whether $|\mu|>1$ or $|\mu|<1$ respectively.

\begin{theorem}
Let $E^n$ denote the ellipse
\[
E^n=\{ \rme^{i\theta}+\rme^{-2ng-i\theta}:\theta\in[-\pi,\pi]\}
\]
and let
\begin{equation}
E_\alp=\{ \lam:  b_{11}(\lam)+b_{22}(\lam)\in E^n
\}.\label{7E}
\end{equation}
Then $B$ has an eigenvalue of modulus $1$ if an only if $\lam\in
E_\alp$. Moreover $E_\alp$ is closed and bounded with
\[
E_\alp\subseteq \Spec(H).
\]
\end{theorem}

\Proof If $\mu_1=\rme^{i\theta}$ for some $\theta\in [-\pi,\pi]$
then $\mu_2=\rme^{-2ng-i\theta}$, and
\[
b_{11}(\lam)+b_{22}(\lam)= \rme^{i\theta}+\rme^{-2ng-i\theta}.
\]
or equivalently $\lam\in E_\alp$. The converse also holds. Our
comments above on the degrees of $b_{ij}(\lam)$ imply that
$|\mu_1+\mu_2|$ increases indefinitely as $|\lam|$ grows. Therefore
one of the $\mu_i$ must have modulus greater than $1$ for large
enough $|\lam|$ and such $\lam$ cannot lie in $E_\alp$; therefore
$E_\alp$ must be bounded. The fact that $E_\alp$ is closed follows
directly from its definition.

Corresponding to any $\lam\in E_\alp$ there exists a solution $f$
of (\ref{import}) such that $f_{m+n}=\rme^{i\theta}f_m$ for some
$\theta\in\R$ and all $m\in\Z$. This $f$ is bounded but its $l^2$
norm is infinite. If we put
\[
f_{\eps,m}=\rme^{-\eps|m|}f_m
\]
then a direct and well-known calculation shows that $\norm f_\eps
\norm_2\to\infty$ and
\[
\frac{\norm (H_0+W)f_\eps -\lam f_\eps \norm_2}{\norm f_\eps
\norm_2}\to 0
\]
as $\eps\to 0$. Applying Theorem 1 we deduce that
\[
\lam\in\Spec(H_0+W)\subseteq \Spec (H).
\]

The set $\C\,\backslash E_\alp$ is the union of disjoint components
and the number of eigenvalues $\mu_j$ of $B$ which have modulus
less than $1$ cannot change within each component, because the
eigenvalues depend continuously on $\lam$. This number must be
either $1$ or $2$, and within the unbounded component it is $1$.
The following theorem joins the components into two sets.

\begin{theorem}
If $\lam$ lies in
\begin{equation}
I_\alp=\{ \lam: b_{11}(\lam)+b_{22}(\lam)\in {\rm int}(E^n )
\}\label{7I}
\end{equation}
then all solutions of (\ref{import}) are exponentially decreasing.
If, however, $\lam$ lies in
\begin{equation}
O_\alp=\{ \lam:  b_{11}(\lam)+b_{22}(\lam)\in {\rm ext}(E^n )
\}\label{7O}
\end{equation}
then there is an exponentially increasing solution of
(\ref{import}). The three sets $I_\alp$, $O_\alp$ and $E_\alp$ are
disjoint and cover $\C$.
\end{theorem}

\Proof The condition (\ref{7I}) holds if and only if both
$\mu_i$ have modulus less than $1$, and this implies that every
solution of (\ref{import}) is exponentially decreasing on $\Z$.
Similarly The condition (\ref{7O}) holds if and only if one $\mu_i$
has modulus greater than $1$, and this implies that one non-zero
solution of (\ref{import}) is exponentially increasing on $\Z$.

The explicit description of the above sets depends upon the value
of $n$. For $n=1$ we have $\alp\in M$ and
\[
E_\alp =E+\alp.
\]
If $n=2$ and $\alp=(\alp_0,\alp_1)\in M^2$ then
\[
B=\left[\begin{array}{cc}-\rme^{-2g}&-\rme^{-3g}(\lam-\alp_1)\\
\rme^{-g}(\lam-\alp_0)&\rme^{-2g}\{(\lam-\alp_0)(\lam-\alp_1)-1\}
\end{array} \right]
\]
and $E_\alp$ is the set of $\lam$ such that
\begin{equation}
\rme^{-2g}\{ (\lam-\alp_0)(\lam-\alp_1)-2\} \in E^2.\label{E2}
\end{equation}
This equation may be solved to present $\lam$ explicitly as a
function of $\theta$. For larger values of $n$ it is probably only
practicable to find $E_\alp$ numerically.

The special case $n=p=1$ of the following theorem was proved in
\cite{EBD}. The idea owes much to the theory of block Toeplitz
matrices \cite{Boe,Boe2,BS,RT}.

\begin{theorem}
Let $H$ be defined by (\ref{schr}) where $g>0$ and $V$ is a
$(\Z,M)$ pseudo-ergodic potential. If $\alp\in M^n$ and $\bet\in
M^p$ then
\[
I_\alp\cap O_\bet \subseteq \sigl (H).
\]
\end{theorem}

\Proof  We consider the operator $K=H_0+W$ acting on $l^2(\Z)$ where
\[
W_m=\left\{ \begin{array}{ll}
\alp_r & \mbox{if $m\geq 0$ and $m=r$ mod $n$}\\
\bet_r & \mbox{if $m<0$ and $m=r$ mod $p$}.
\end{array} \right.
\]
We then consider the solutions of
\[
\rme^{-g} f_{m-1} +W_m f_m +\rme^g f_{m+1}=\lam f_m
\]
where $\lam \in I_\alp \cap O_\bet$. Since $\lam\in O_\bet$ there
exists a solution $f$ which is exponentially growing for $m<0$,
i.e. which decreases exponentially as $m\to -\infty$. Continuing
this solution to positive $m$ it follows from $\lam\in I_\alp$ that
$f$ also decreases exponentially as $m\to\infty$. Hence $f$ is an
eigenvector of finite variance and $\lam\in \sigl(K)\subseteq
\sigl(H)$.

\begin{theorem}
If in addition to the hypotheses of the last theorem we put $M=
[-\mu,\mu]$ then $\Spec(H)=E+[-\mu,\mu]$ for all $\mu\geq
\rme^g+\rme^{-g}$. Moreover $0\in
\Spec(H)$ if and only if $\mu
\geq\rme^g-\rme^{-g}$.
\end{theorem}

\Proof The first statement only  needs the observation that the two
sides of (\ref{nn}) coincide under the given condition. If $\mu <
\rme^g-\rme^{-g}$ then $0\notin\Spec(H)$ by Theorem 4. Now $0\in
E_{(-\mu,\mu)}$ if and only if $\rme^{-2g}(-\mu^2-2)\in E^2$ by
(\ref{E2}), and this is equivalent to $\mu=\rme^g-\rme^{-g}$; for
such $\mu$ one has $0\in\Spec(H)$ by Theorem 10. For smaller $\mu$
we have $0\in I_{(-\mu,\mu)}$ and for larger $\mu$ we have $0\in
O_{(-\mu,\mu)}$. Therefore $0\in I_{(0)}\cap O_{(-\mu,\mu)}$ for
$\mu>\rme^g-\rme^{-g}$, and $0\in\sigl(H)$ by Theorem 12.

If $M=[-\mu,\mu]$ the above theorems admit the possibility that
there are two holes in the spectrum on either side of the origin
for
\[
\rme^g-\rme^{-g}<\mu< \rme^g+\rme^{-g}.
\]
We nevertheless conjecture that one has $\Spec(H)=\conv(E+M)$ for
all $\mu\geq\rme^g-\rme^{-g}$.

We contrast the above with the case in which $M=\{\pm\mu\}$. The
following theorem completely determines the real part of $\Spec(H)$
under the stated conditions.

\begin{theorem}
If $M=\{\pm \mu\}$ and $\mu >\rme^g+\rme^{-g}$ then
\[
(\conv(E)+\mu)\cup(\conv(E)-\mu)\subseteq \Spec(H) \subseteq
B(\mu,\rme^{g}+\rme^{-g})\cup B(-\mu,\rme^{g}+\rme^{-g})
\]
and
\[
\Spec(H) \subseteq \conv(E)+[-\mu,\mu].
\]
\end{theorem}

\Proof The first inclusion of the statement follows from the case
$n=p=1$ of Theorem 12 as in \cite{EBD}. The second follows from the
first half of Theorem 4, and the final one follows from Theorem 3.

We conjecture that the first inclusion is actually an equality.

\begin{corollary}
If $M=\{\pm \mu\}$ then $0\in\Spec (H)$ if and only if
\[
\rme^g-\rme^{-g}\leq \mu\leq \rme^g+\rme^{-g}.
\]
\end{corollary}

\Proof If $\mu< \rme^g-\rme^{-g}$ then $0\notin\Spec(H)$ by
combining Corollary 2 and Theorem 13. If $\mu > \rme^g+\rme^{-g}$
then $0\notin\Spec(H)$ by Theorem 14. If $\mu= \rme^g-\rme^{-g}$
then $0\in E_{(-\mu,\mu)}\subseteq \Spec(H)$ by Theorem 10. If
$\mu=\rme^g+\rme^{-g}$ then $0\in E_{\mu}\subseteq \Spec(H)$ by
Theorem 10. Finally if $\rme^g-\rme^{-g}< \mu <
\rme^g-\rme^{-g}$ then $0\in O_{(-\mu,\mu)}\cap I_\mu \subseteq
\Spec (H)$ by Theorem 12.

We next turn to the nsa Anderson model in $\Z^n$. The operator $H$
on $l^2(\Z^n)$ is defined by
\[
(Hf)(m,n)=(H_0f)(m,n)+V(m,n)f(m,n)
\]
where
\[
(H_0f)(m,n)=\rme^g f(m+1,n)+\rme^{-g}f(m-1,n)+f(m,n+1)+f(m,n-1)
\]
for some $g>0$. We assume that $V$ is real-valued and
pseudo-ergodic with values in $M=[-\mu,\mu]$. It follows by Fourier
transform methods that $H_0$ is normal with spectrum equal to
\[
\tilde{E}=E+[-2(n-1),2(n-1)]
\]
where $E$ is the ellipse defined by (\ref{ellipse}). This set is
connected with a hole around the origin if $n=2$ but it may or may
not have such a hole for $n\geq 3$. This phenomenon is a result of
the particular choice of lattice used to discretize the Laplacian.
If $\mu$ is sufficiently small the same applies to $\Spec(H)$.

\begin{theorem}
If $\mu\geq \rme^g+\rme^{-g}-2(n-1)$ then $\Spec(H)$ is the convex
set
\[
E+[-\mu-2(n-1),\mu+2(n-1)].
\]
\end{theorem}

\Proof As in the one-dimensional case we need only observe that the two sides of
(\ref{nn}) are equal under the hypotheses.

We next mention the same operator acting in $l^2(X)$ where
\[
X=\{(m,n):m\in \Z,1\leq n\leq N\}
\]
subject to Dirichlet boundary conditions; the Neumann case is
similar. We may carry out an analysis similar to that above if we
are only concerned to determine the spectrum, but more detailed
spectral information is obtained by putting $l^2(X)=l^2(\Z,\cK)$
where $\cK=\C^N$. We then put
\[
(H_0f)(m)=\rme^g f(m+1)+\rme^{-g}f(m-1)
\]
and
\[
\tilde{V}(m)(r,s)=\left\{
\begin{array}{ll}
1&\mbox{if $|r-s|=1$}\\ V(m,r)&\mbox{if $r=s$}\\0&\mbox{otherwise}
\end{array}
\right.
\]
where $1\leq r,s \leq N$ in all cases. Note that $H_0$ is normal
and $\tilde{V}(m)$ is a self-adjoint matrix for all $m\in
\Z$, so all of the theorems of Section 2 apply. Using such ideas it
is possible to analyze the localized spectrum of $H$ as in the
one-dimensional case.

We finally comment that certain random {\it bidiagonal} operators
can also be treated by the methods of this paper by making the
appropriate choice of $H_0$, as can a variety of other operators
whose matrix coefficients depend only on $m-n$ whenever $m\not= n$.
See \cite{BZ,FZ,TCE}, which use probabilistic rather than
pseudo-ergodic methods.

\section{Resolvent Norms}

The spectral behaviour of a bounded operator $A$ acting on a
Hilbert space $\cH$ can be measured in several ways. In
pseudospectral theory one examines the contours of the function
\[
s(A,z)=\left\{ \begin{array}{ll}
\norm (A-z)^{-1}\norm^{-1}&\mbox{if $z\notin \Spec (A)$}\\
0&\mbox{if $z\in\Spec (A)$.}
\end{array} \right.
\]
This function converges to zero as $z$ approaches the spectrum of
$A$ because of the upper bound
\[
s(A,z)\leq \dist(z,\Spec(A))
\]
and the case of most interest is when $s(A,z)$ is very small for
$z$ far from the spectrum. The determination of the pseudospectra,
defined as the family of sets $\{  z: s(z)<\eps \}$ for all
positive $\eps$, is computationally heavy, but the family carries
much more information than the spectrum alone \cite{Boe, BS,
Red,RT,Tre1,Tre2}.

\begin{lemma}
The function $s(A,\,\cdot\,)$ satisfies the Lipschitz inequality
\[
|s(A,z)-s(A,w)|\leq |z-w|
\]
for all $z,w\in \C$.
\end{lemma}

The proof uses the formula
\begin{equation}
s(A,z)=\inf\{\norm (A-z)f\norm/\norm f\norm :0\not=
f\in\cH\}\label{sform}
\end{equation}
valid for all $z\notin\Spec(A)$. Note that this may be false for
$z\in\Spec (A)$, as one may see by considering the operator
$\hat{A}$ on $l^2(\Z^+)$ defined by
\[
\hat{A}f(n)=\left\{ \begin{array}{ll}
 0&\mbox{if $n=1$}\\
f(n-1)&\mbox{if $n\geq 2$.}
\end{array} \right.
\]
The next theorem provides an upper bound on $s(A,\,\cdot\,)$ which
may be used to compute it numerically. Let $L$ be a
finite-dimensional subspace of $\cH$ and let $P$ be the orthogonal
projection onto $L$. We define $B(A,L,z)$ to be the restriction of
\[
P(A-zI)^\ast P(A-zI)P +PA^\ast (I-P)AP
\]
to the subspace $L$, and $\sig(A,L,z)$ to be the square root of the
smallest eigenvalue of $B(A,L,z)$.

\begin{theorem}
If we put
\[
s(A,L,z)=\min\{ \sig(A,L,z),\sig(A^\ast,L,\overline{z})\}
\]
then
\[
s(A,L,z)\geq s(A,z).
\]
The functions $s(A,L,\,\cdot\,)$ decrease monotonically and locally
uniformly to $s(A,\,\cdot\,)$ as the subspaces increase.
\end{theorem}

\Proof It follows from its definition that
\[
\sig(A,L,z)=\min\{ \norm (A-z)f\norm/\norm f\norm :0\not= f\in L\}.
\]
It is clear from this that $\sig (A,L,z)$ decreases monotonically
and pointwise to
\[
\sig(A,z)=\inf\{ \norm (A-z)f\norm/\norm f\norm :0\not= f\in \cH\}.
\]
If $z\notin\Spec (A)$ this equals $s(A,z)$. Similar comments apply
with $A$ replaced by $A^\ast$, and we also have
$\sig(A,z)=\sig(A^\ast,\overline{z})$ for all $z$.

On the other hand if $z\in\Spec (A)$ we have either $\sig(A,z)=0$
or $\sig(A^\ast,\overline{z})=0$, or both. This implies that
$s(A,L,z)$ converges monotonically and pointwise to $0$. Since all
the functions involved are Lipschitz continuous with Lipschitz
constant $1$, the convergence must be  locally uniform.

Now suppose that $\cH$ equals $l^2(X,\cK)$ and $L$ is defined as
the space of all functions with support in a particular finite
region $\Ome$. The above theorem is better than the mere
computation of the spectrum of $PAP$ restricted to $L$ (possibly
subject to certain boundary conditions on $\partial\Ome$) because
it gives rigorous upper bounds to $s(A,z)$ rather than uncontrolled
approximations. Another advantage is that it provides an upper
bound for $s(A,z)$ for {\it every} extension of the operator $A$
beyond the subspace $L$. Because of its approximate nature one
cannot determine the spectrum of $A$ exactly using the above
theorem, but it may be possible to get good approximations to the
pseudospectra, which are often of greater importance for such
operators.

We now turn to pseudo-ergodic operators, working in the technical
context of Section 2. The following theorem indicates how one may
get rigorous upper bounds and approximations to the pseudospectra
by selecting appropriate potentials $W$.

\begin{theorem} If $H=H_0+V$ where $V$ is $(\Gam,M)$ pseudo-ergodic
and $K=H_0+W$ where $W:X\to M$ is arbitrary, then
\[
\norm (H-zI)^{-1}\norm \geq \norm (K-zI)^{-1}\norm
\]
for all $z\in \C$. Therefore
\[
s(H,z)=\min\{ s(H_0+W,z):W\in M^X\}.
\]
If $V,W$ are both $(\Gam,M)$ pseudo-ergodic then the resolvent
norms and hence pseudospectra of $H$ and $K$ are equal.
\end{theorem}

\Proof By Theorem 1 we need only consider the case in which $z$
does not lie in the spectrum of either operator. If $s(K,z)<c$ then
there exists $f\in l^2(X,\cK)$ such that $\norm (K-z)f\norm<c\norm
f\norm$ and by approximation we may assume that $f$ has finite
support. Using the pseudo-ergodic property of $H$ there exists $g$
of finite support such that $\norm (H-z)g\norm<c\norm g\norm$ and
this implies that $s(H,z)<c$. Hence $s(H,z)\leq s(K,z)$. The
remainder of the proof follows Theorem 1 or Corollary 2.

For the nsa periodic Anderson model with $M=[-\mu,\mu]$ the
asymptotic limit of the finite volume spectrum has been determined
\cite{Gold}, and it is seen that for certain ranges of the
parameter $\mu$ zero does not lie in the asymptotic spectrum, which
is the union of a set of complex curves. On the other hand the
spectrum of the same operator on any finite interval subject to
Dirichlet boundary conditions is entirely real. It has been
suggested in \cite{Dahmen} that for periodic boundary conditions
there is no pseudospectral pathology of the type which occurs for
Dirichlet boundary conditions. However, our results demonstrate
that spatially rare special sections of a random potential have a
dominant effect on the spectrum of the infinite volume nsa Anderson
operator. This should not be taken as an indication that our
results are unphysical: it is well known that the behaviour of bulk
materials is often radically affected by the presence of low
concentrations of impurities and/or defects, and one should expect
the mathematics to reflect this.

We have implemented the above ideas numerically using Matlab for
the operator $H$ defined by (\ref{schr}) where $\rme^g=2$ and $V_n$
are independent random variables uniformly distributed on $[-3,3]$.
We took $L$ to be the subspace of all sequences with support in
$[1,100]$ and computed the minimum value of $\sig (H,L,x)$ over
$1000$ different choices of the potential $V$. We chose to study
real $x\in [0,6]$, but complex values of $x$ in any region can be
accommodated by the same method. This yielded the upper bounds
$\overline{s}(H,x)$ as follows (the omitted values of
$\overline{s}(H,x)$ all vanish to the given accuracy).

\[
\begin{array}{cc}
x&\overline{s}(H,x)\\

\hline

0.0  & 0.0283 \\

0.5 &  0.0203 \\

1.0  &  0.0084\\

1.5  &  0.0015\\

...&...\\

4.5  &0.0044\\

5.0& 0.2233\\

5.5&0.6259\\

6.0&1.0817

\end{array}
\]

Our general theory shows that the real part of the spectrum of this
operator is $[-5.5,5.5]$, which is consistent with the numerical
conclusion that
\[
\norm (H- xI)\inv\norm\geq 10^{2}
\]
for all $1.0\leq x\leq 4.5$ and
\[
\norm (H- xI)\inv\norm\geq 10^{4}
\]
for all $2.0\leq x\leq 4.0$. (Of course the numerical calculation
can also be carried out in cases in which one does not have a prior
theoretical solution!) The eigenvectors of $B(A,L,x)$ corresponding
to the smallest eigenvalues were also computed for several values
of $x$. As expected from the theory of localized spectrum, they
were all highly concentrated around some point in the interior of
$[1,100]$, and negligible at the ends of the interval.

We finally examine the behaviour of the resolvent norm at the point
$z=0$. To be precise we consider the nsa Anderson model with
$M=[-\mu,\mu]$ acting on $l^2(\Z)$ for various values of $\mu$.
Recall that Theorem 13 states that $0\in\Spec (H)$ if and only if
$\mu\geq\rme^g -\rme^{-g}$.

\begin{theorem}
If  $\,\,0\leq \mu < \rme^g -\rme^{-g}$ then
\[
\norm H^{-1}\norm^{-1}= \rme^g -\rme^{-g}-\lam.
\]
\end{theorem}

\Proof If we exhibit the $\mu$ dependence of $H$ explicitly and
put $ t(\mu)=s(H_\mu,0)$ then it follows from (\ref{sform}) that
\[
|t(\mu)-t(\nu)|\leq |\mu-\nu|
\]
for any $0\leq \mu,\nu <\rme^g -\rme^{-g}$. Since $t(0)=\rme^g
-\rme^{-g}$ and $t(\rme^g -\rme^{-g})=0$ we conclude that
\[
t(\mu)=\rme^g -\rme^{-g}-\mu
\]
for all $0\leq \mu <\rme^g -\rme^{-g}$.

\section{Constrained Potentials}

We have avoided the use of any probabilistic methods by the
introduction of the concept of pseudo-ergodicity. We now explore
the variety of situations in which our ideas are applicable. The
obvious possibility is to assume that $\mu$ is a probability
measure with support equal to the set $\cM\subseteq \cL(\cK)$ and
to assume that $V_x$ are independent random variables as $x\in X$
varies and that each is distributed according to $\mu$. However,
even if we assume that $V_x$ are independent, we may permit each
$V_x$ to be distributed according to a different probability
measure $\mu_x$ with support equal to $\cM$. These measures need
not even be $\Gam$-stationary, but they must satisfy the following
condition. For every open set $U\subset \cL(\cK)$ such that $U\cap
\cM\not=
\emptyset$ there must exist a constant $c_U>0$ such that $\mu_x(U)\geq
c_U$ for all $x\in X$. This is sufficient to imply that $V$ is
$(\Gam,\cM)$ pseudo-ergodic almost surely by the usual
probabilistic argument. For all such probabilistic models the
spectrum (or localized spectrum) of the operator $H$ is the same.

Similar remarks apply to a variety of other probabilistic models in
which the values $V_x$ are not independent. There is one situation,
however, in which changes in the spectrum may arise. We say that a
potential $V$ satisfies the local constraints
$\cQ=(\cM,\gam_1,...,\gam_k,\cN_1,...,\cN_k)$ where $\gam_i\in\Gam$
and $\cM,\cN_i$ are closed, bounded subsets of $\cL(\cK)$ under the
following conditions. For all $x\in X$ we require that $V_x \in
\cM$ and also that
\[
V_x-V_{\gam_i x}\in \cN_i
\]
for all $i=1,...,k$. Even more general constraints can be
formulated. We then say that $V$ is $(\Gam,\cQ)$ pseudo-ergodic if
it satisfies the constraints $\cQ$ and for any other potential $W$
which satisfies the same constraints and any finite subset $F$ of
$X$ and any $\eps >0$ there exists $\gam\in\Gam$ such that
\[
\norm V_{\gam x}-W_x \norm <\eps
\]
for all $x\in F$. These constraints force a relationship between
the values of $V_x$ at neighbouring points which is stronger than a
mere probabilistic correlation.

\begin{lemma} If $H_j=H_0+V_j$ where $V_1$ is $(\Gam,\cM)$ pseudo-ergodic
and $H_2$ is $(\Gam,\cQ)$ pseudo-ergodic, then
\[
\Spec(H_2)\subseteq \Spec(H_1).
\]
Any two $(\Gam,\cQ)$ pseudo-ergodic operators have the same
spectrum.
\end{lemma}

\Proof The first statement is a consequence of Theorem 1. The
second involves adapting the proof of the same theorem.

We now apply the above ideas in a simple context. We assume that
$X=\Z$, that $\Gam$ is the usual translation group acting on $\Z$,
and that $\cK=\C$. We assume that $a,b$ are two positive constants
and impose attractive constraints $\cQ_1$ of the form
\[
-a\leq V_n\leq a \hskip 0.3in , \hskip 0.3in |V_n-V_{n+1}|\leq b
\]
for all $n\in\Z$. Although we are not able to prove Theorem 12 in
full generality under such conditions the important special case
$n=p=1$ is still valid.

\begin{theorem}
Let $H$ be defined by (\ref{schr}) where $g>0$ and $V$ is a
$(\Z,\cQ_1)$ pseudo-ergodic potential. We have
\[
E+[-a,a]\subseteq \Spec (H) \subseteq \conv(E)+[-a,a]
\]
where $E$ is given by (\ref{ellipse}). If $\alp,\bet\in [-a,a]$
then
\[
I_\alp\cap O_\bet \subseteq \sigl (H).
\]
\end{theorem}

\Proof  The first statement of the theorem is proved as in Theorem 3.
For the second part we follow the method of Theorem 12 but for the
operator $K=H_0+W$ acting on $l^2(\Z)$ where
\[
W_n=\left\{ \begin{array}{ll}
\alp & \mbox{if $n> N$}\\
\bet & \mbox{if $n< 0$}\\
\bet+n(\alp-\bet)/N & \mbox{if $0\leq n\leq N$.}
\end{array} \right.
\]
Here we take $N$ large enough to ensure that $W$ satisfies the
constraints $\cQ_1$.

A more interesting variation upon our earlier theory occurs if we
impose the repulsive constraint $\cQ_2$ defined by
\[
-a\leq V_n\leq a \hskip 0.3in , \hskip 0.3in |V_n-V_{n+1}|\geq b
\]
for all $n\in\Z$, where $0<b\leq 2a$. This excludes constant
potentials, thus rendering the first inclusion of Theorem 3
invalid. The range of a $(\Gam,\cQ_2)$ pseudo-ergodic potential $V$
is equal to $M=[-a,a-b]\cup [b-a,a]$.

The spectrum of the Anderson model (\ref{schr}) is easy to
determine in the self-adjoint case, and we start with this.

\begin{theorem}
If $g=0$ and $V$ is $(\Gam,\cQ_2)$ pseudo-ergodic then the spectrum
of the operator $H=H_0+V$ defined by (\ref{schr}) is given by
\[
\Spec (H)=T\cup(-T)
\]
where
\[
T=\left[ b-a,a-\frac{b}{2}+\sqrt{\frac{b^2}{4}+4}\right].
\]
Thus $\Spec(H)$ has a spectral gap if and only if $a<b\leq 2a$.
\end{theorem}

\Proof Let $W$ be the potential $W_n=(-1)^nb/2$, so that $W+sI$
satisfies the constraints $\cQ_2$ for all real $s$ such that
$|s|\leq a-b/2$. It follows from Theorem 1 that if $S=\Spec
(H_0+W)$ then
\begin{equation}
S+[b/2-a,a-b/2]\subseteq \Spec (H).\label{impo}
\end{equation}
Conversely $\norm V-W\norm \leq a-b/2$, so the perturbation
theoretic argument used in Theorem 4 implies that
\[
\Spec(H)\subseteq S+[b/2-a,a-b/2].
\]
We deduce that
\[
\Spec(H)= S+[b/2-a,a-b/2]
\]
and complete the proof by using a Bloch wave analysis to compute
the set $S$.

Now let us denote the same operator by $L_g$ for $g \geq 0$. We may
regard $L_g$ as a perturbation of $L_0$ and use the argument of
Theorem 4 to show that
\[
\Spec(L_g)\subseteq \Spec (L_0) +B(0,\rme^g - 1).
\]
We may also use Theorem 4 as it stands to obtain an outer estimate
of $\Spec (L_g)$. We may obtain inner estimates by the method of
Section 3 provided we are careful to avoid the use of constant
potentials.

\begin{theorem}
We have
\[
S+\left[ \frac{b}{2}-a,a-\frac{b}{2}\right]\subseteq \Spec (L_g)
\]
where
\[
S=\left\{
\pm\sqrt{   \frac{b^2}{4}+2+\rme^{2g+i\theta}+\rme^{-2g-i\theta} }
:\theta \in [-\pi,\pi]
\right\}.
\]
\end{theorem}

\Proof if we put $\alp_0=-b/2$ and $\alp_1=b/2$ and solve (\ref{E2}) for $\lam$
we obtain
\[
E_{(-b/2,b/2)}=S.
\]
The remainder of the proof follows Theorem 23, using the last part
of Theorem 10.

Note that for small positive $g$, $S$ consists of two closed curves
on opposite sides of the $y$-axis, but for large $g$ it is a single
curve enclosing the origin.

\vskip 1in
{\bf Acknowledgements} I acknowledge valuable conversations with N
Trefethen and I Goldsheid during the course of this work. I also
thank the EPSRC for support under grant no. GR/L75443.

\vskip 1in

Department of Mathematics\newline%
King's College\newline%
Strand\newline%
London WC2R 2LS\newline%
England
\vskip 0.2in
e-mail: E.Brian.Davies@kcl.ac.uk
\vfil
\end{document}